\documentclass[aop,MSNbibl,seceqn,dvips]{arximspdf}

\doi{10.1214/12-AOP800} 
\volume{42}
\issue{3}
\pubyear{2014}
\firstpage{946}
\lastpage{958}

\makeatletter
\def\cal{\mathcal}
\def\eps{\varepsilon}
\def\P{{\cal P}}
\def\LL{{\cal L}}
\def\FF{{\cal F}}
\def\eps{{\varepsilon}}
\def\ch{\operatorname{ch}}
\def\vec{}
\newcommand{\e}{\mathbb{E}}
\newcommand{\Reals}{\mathbb{R}}
\newcommand{\Natural}{\mathbb{N}}
\newcommand{\vsi}{{\vec{\sigma}}}
\newcommand{\vrho}{{\vec{\rho}}}

\newtheorem{proposition}{Proposition}
\newtheorem{lemma}{Lemma}
\newtheorem{theorem}{Theorem}

\makeatother

\begin{document}
\begin{frontmatter}

\title{The Parisi formula for mixed \lowercase{$p$}-spin models}
\runtitle{Parisi formula}

\begin{aug}
\author[A]{\fnms{Dmitry} \snm{Panchenko}\corref{}\ead[label=e1]{panchenk@math.tamu.edu}\thanksref{t1}}
\thankstext{t1}{Supported in part by NSF grant.}
\runauthor{D. Panchenko}
\affiliation{Texas A\&M University}
\address[A]{Department of Mathematics\\
Texas A\&M University\\
Mailstop 3386\\
College Station, Texas 77843\\
USA\\
\printead{e1}}
\end{aug}

\received{\smonth{2} \syear{2012}}
\revised{\smonth{9} \syear{2012}}

%
\begin{abstract}
The Parisi formula for the free energy in the Sherrington--Kirk\-patrick
and mixed $p$-spin models for even $p\geq2$ was proved in the seminal
work of Michel Talagrand [\textit{Ann. of Math.} (2) \textbf{163}
(2006) 221--263].
In this paper we prove the Parisi
formula for general mixed $p$-spin models which also include $p$-spin
interactions for odd $p$. Most of the ideas used in the paper are well
known and can now be combined following a recent proof of the Parisi
ultrametricity conjecture in [\textit{Ann. of Math.} (2) \textbf{177}
(2013) 383--393].
\end{abstract}

%
\begin{keyword}[class=AMS]
\kwd{60K35}
\kwd{82B44}
\end{keyword}
\begin{keyword}
\kwd{Sherrington--Kirkpatrick model}
\kwd{free energy}
\kwd{ultrametricity}
\end{keyword}

\end{frontmatter}

\section{Introduction and main result}
The formula for the free energy in the Sherrington--Kirkpatrick model
\cite{SK} was famously discovered by G. Parisi in~\cite{Parisi79,Parisi}
using the approach that combined a replica trick with a very
special choice of the replica matrix. It was later understood in \cite
{M1,M2} that the special form of the replica matrix conjectured
by Parisi corresponded to a number of physical properties of the Gibbs
measure of the model, one of them being the ultrametricity of its
support. The Parisi formula for the free energy in the
Sherrington--Kirkpatrick and mixed $p$-spin models was proved by M.
Talagrand in \cite{T-P} following the discovery of the replica symmetry
breaking interpolation scheme by F. Guerra in \cite{Guerra}. However,
for technical reasons only the case of $p$-spin interactions for even
$p\geq2$ was considered. Using the main result in \cite{PU}, which
yields that under a small perturbation of the Hamiltonian the support
of the Gibbs measure in these models is indeed asymptotically
ultrametric, we prove the Parisi formula for general mixed $p$-spin
models that include odd $p$-spin interactions as well.

Let $N\geq1$. Let us consider Gaussian processes $H_{N,p}(\vsi)$ for
$p\geq1$ indexed by $\vsi\in\Sigma_N = \{-1,+1\}^N$, called pure
$p$-spin Hamiltonians,
%
%
\begin{equation}
H_{N,p}(\vsi) = \frac{1}{N^{(p-1)/2}} \sum_{1\leq i_1,\ldots
,i_p\leq N}g_{i_1,\ldots,i_p}
\sigma_{i_1}\cdots\sigma_{i_p}, \label{HamSKp}
\end{equation}
where random variables $(g_{i_1,\ldots,i_p})$ are standard Gaussian
independent for all \mbox{$p\geq1$} and all $(i_1,\ldots,i_p)$. Let us define
a mixed $p$-spin Hamiltonian as their linear combination
%
%
\begin{equation}
H_N(\vsi) = \sum_{p\geq1}
\beta_p H_{N,p}(\vsi) \label{HN}
\end{equation}
with coefficients $(\beta_p)$ that decrease fast enough, for example,
$\sum_{p\geq1} 2^p \beta_p^2<\infty$. This technical condition is
sufficient to ensure that the process is well defined when the sum
includes infinitely many terms. The covariance of the Gaussian process
$H_N(\vsi)$ is easy to compute and is given by a function of the
normalized scalar product, called overlap, $R_{1,2}= N^{-1}\sum_{i\leq
N} \sigma^1_i \sigma_i^{2}$ of spin configurations $\vsi^1$ and~$\vsi^2$,
%
%
\begin{equation}
\e H_N \bigl(\vsi^1 \bigr) H_N \bigl(
\vsi^2 \bigr) = N\xi(R_{1,2}),
\end{equation}
where $\xi(x)=\sum_{p\geq1}\beta_p^2 x^p$. Given $k\geq1,$ let us
consider two sequences of parameters,
%
%
\begin{equation}
0\leq m_0\leq m_1\leq\cdots\leq m_{k-1}\leq
m_k \leq1 \label{m}
\end{equation}
and
%
%
\begin{equation}
0=q_0\leq q_1\leq\cdots\leq q_{k}\leq
q_{k+1} = 1, \label{q}
\end{equation}
which will be denoted by $\vec{m}$ and $\vec{q}$, and consider
independent Gaussian random variables $(z_j)_{0\leq j\leq k}$ with
variances $\e z_j^2 = \xi'(q_{j+1}) - \xi'(q_j).$
We define
%
%
\begin{equation}
X_{k+1}=\log\ch\sum_{0\leq j\leq k} z_j
\quad\mbox{and}\quad X_l=\frac{1}{m_l}\log\e_l\exp
m_l X_{l+1} \label{Xl}
\end{equation}
recursively for $l\leq k,$ where $\e_l$ denotes the expectation in the
r.v. $(z_j)_{j\geq l}.$
When $m_l=0$ this means that $X_l=\e_l X_{l+1}.$ Let us denote $\theta
(q)=q\xi'(q)-\xi(q)$ and define
%
%
\begin{equation}
\P_k(\vec{m},\vec{q}) = \log2 + X_0(\vec{m},\vec{q}) -
\frac{1}{2} \sum_{1\leq j\leq k} m_j \bigl(
\theta(q_{j+1}) - \theta(q_j) \bigr). \label{Pk}
\end{equation}
Then the following theorem holds.
%
%
\begin{theorem}[(The Parisi formula)]\label{Th1} We have
%
%
\begin{equation}
\lim_{N\to\infty} \frac{1}{N} \e\log\sum_{\vsi\in\Sigma_N}
\exp H_N(\vsi) = \inf\P_k(\vec{m},\vec{q}),
\label{Parisi}
\end{equation}
where the infimum is taken over all $k, \vec{m}$ and $\vec{q}$ as above.
\end{theorem}

The quantity in the limit on the left-hand side is called the free
energy of the model and the infimum on the right-hand side is the
famous Parisi formula. One can include the external field term in the
model, but for simplicity of notation we will omit it. The proof we
give here, obviously, assumes a certain level of expertise, but all the
details starting from the foundations can be found in \cite{SKmodel}.

\section{Proof}

Most of the ideas of the proof are well known and available in
different places in the literature. Under various formulations of the
ultrametricity conjecture, one can find arguments that contain many of
the same ideas in \cite{ACh} and \cite{Pspins} in the case of models
with only even $p$-spin interactions, and a sketch of the proof of the
general case in Section 15.3 in \cite{SG2}. The ingredient that was
missing is the main result in \cite{PU} which also allows us to handle
the case of the general mixed $p$-spin models.

\textit{The Ghirlanda--Guerra identities}. A central role in the proof
is played by the Ghirlanda--Guerra identities \cite{GG} that are
utilized in two distinct ways. First, they yield positivity of the
overlap via Talagrand's positivity principle, which allows us to obtain
the upper bound using Guerra's replica symmetry breaking interpolation
scheme and, second, they imply ultrametricity of the overlap array
using the main result in \cite{PU}, which allows us to identify the
asymptotic Gibbs measures that appear in the proof of the lower bound
based on the Aizenman--Sims--Starr scheme~\cite{AS2}. Let us consider a
perturbation Hamiltonian
%
%
\begin{equation}
H_N^{\mathrm{pert}}(\vsi) = N^{-1/8} \sum
_{p\geq1} 2^{-p} x_p H_{N,p}'(
\vsi), \label{Hpert}
\end{equation}
where $H_{N,p}'(\vsi)$ are independent copies of the $p$-spin
Hamiltonians in (\ref{HamSKp}) and $(x_p)_{p\geq1}$ are i.i.d. random
variables uniform on an interval of length one, for example, $[1,2]$.
Replacing $H_N$ with $H_N + H_N^{\mathrm{pert}}$ in (\ref{Parisi}),
obviously, does not affect the limit since the perturbation term is of
a smaller order. However, adding this perturbation term regularizes the
Gibbs measure in the following way. Let $G_N$ be the Gibbs measure on
$\Sigma_N$ corresponding to the Hamiltonian $H_N + H_N^{\mathrm{pert}}$,
%
%
\begin{equation}
G_N(\vsi) = \frac{\exp(H_N(\vsi) + H_N^{\mathrm{pert}}(\vsi))}{Z_N},
\end{equation}
where $Z_N = \sum_{\vsi\in\Sigma_N} \exp(H_N(\vsi)+H_N^{\mathrm
{pert}}(\vsi
))$, and denote by $\langle \cdot\rangle $ the average with respect
to the
product Gibbs measure $G_N^{\otimes\infty}$. Let $(\vsi^l)_{l\geq1}$
be an i.i.d. sequence of replicas sampled from $G_N$ and denote by
%
%
\begin{equation}
R_{l,l'}=\frac{1}{N}\sum_{i\leq N}
\sigma^l_i \sigma_i^{l'}
\label{overlap}
\end{equation}
the normalized scalar product, or overlap, of $\vsi^l$ and $\vsi^{l'}.$
Given $p\geq1$, $n\geq2$ and a bounded measurable function $f$ of the
overlaps $(R_{l,l'})_{l,l'\leq n}$ on $n$ replicas, let
%
%
\begin{equation}
\phi(f,n,p) = \Biggl| \e_g \bigl\langle f R_{1,n+1}^p
\bigr\rangle - \frac{1}{n}\e_g \langle f \rangle
\e_g\bigl\langle R_{1,2}^p\bigr\rangle -
\frac{1}{n}\sum_{l=2}^{n}
\e_g \bigl\langle f R_{1,l}^p\bigr\rangle \Biggr|,
\label{GG}
\end{equation}
where $\e_g$ denotes the expectation with respect to all Gaussian
random variables for a fixed uniform sequence $(x_p)_{p\geq1}$.
Then, the Ghirlanda--Guerra identities can be stated as follows.
%
%
\begin{proposition}\label{PropGG}
For any $p\geq1, n\geq2$ and a bounded function $f$ of the overlaps
$(R_{l,l'})_{l,l'\leq n}$,
%
%
\begin{equation}
\lim_{N\to\infty} \e_x \phi(f,n,p) = 0, \label{GGxlim}
\end{equation}
where $\e_x$ is the expectation with respect to $(x_p)_{p\geq1}$.
\end{proposition}

The proof of this result is well known and we refer to Chapter 12 in
\cite{SG2} for details. We will not be using these identities directly
for the measure $G_N$, but for other Gibbs measures with a slightly
modified Hamiltonian $H_{N}(\vsi)$, since it is well known that the
proof of the identities is robust to such modifications and depends
mostly on the form of the perturbation Hamiltonian (\ref{Hpert}). It is
interesting to note that once we finish the proof of Theorem \ref{Th1},
the argument in \cite{PGGmixed} will immediately imply that (\ref
{GGxlim}) holds in a strong sense without the perturbation Hamiltonian
for all $p\geq1$ such that $\beta_p\not=0$ in (\ref{HN}).

\textit{Guerra's replica symmetry breaking bound}.
In the case when $p$-spin interactions for odd $p\geq3$ are not
present in (\ref{HN}), the inequality $\leq$ in (\ref{Parisi}) was
proved by F. Guerra in \cite{Guerra} by inventing the replica symmetry
breaking interpolation scheme. The fact that this inequality holds even
in the presence of odd $p$-spin interactions was observed by M.
Talagrand in \cite{TG} and we will only briefly recall the main idea,
which is to write down Guerra's interpolation scheme in terms of the
Ruelle probability cascades \cite{Ruelle} (Poisson--Dirichlet cascades
in the terminology of \cite{SG2}) and force the overlap to be positive
along the interpolation by adding the perturbation term (\ref{Hpert}).
Given $k\geq1$, the Ruelle probability cascades are defined as (i) a
random probability measure $(w_\alpha)_{\alpha\in\Natural^k}$ on
$\Natural^k$ via some explicit construction involving Poisson processes
on $(0,\infty)$ with the mean measures $\zeta x^{-1-\zeta}\,dx$ for
$\zeta
\in(0,1)$ and (ii) a Gaussian process $(z_\alpha)_{\alpha\in
\Natural
^k}$ with the covariance $\e z_{\alpha^1} z_{\alpha^2} = \xi'(q_{\alpha
^1\wedge\alpha^2})$ where
\begin{eqnarray*}
\alpha^1\wedge\alpha^2 &=& \min \bigl\{l\geq1\dvtx
\alpha_l^1\not= \alpha_l^2 \bigr
\} \qquad\mbox{if } \alpha^1\not=\alpha^2\quad \mbox{and}\\\
\alpha^1
\wedge\alpha^2 &= &k+1 \qquad\mbox{if }\alpha^1=\alpha^2
\end{eqnarray*}
(see Chapter 14 in \cite{SG2} for details). For $0\leq t\leq1$ we
define an interpolating Hamiltonian
%
%
\begin{equation}
H_{N,t}(\vsi,\alpha) = \sqrt{t} H_N(\vsi) + \sqrt{1-t} \sum
_{i\leq N} z_{\alpha, i} \sigma_i,
\label{Hta}
\end{equation}
where $(z_{\alpha,i})_{\alpha\in\Natural^k}$ are independent copies of
$(z_\alpha)_{\alpha\in\Natural^k}$ for $i\geq1$, and let
%
%
\begin{equation}
\varphi(t)=\frac{1}{N} \e\log\sum_{\alpha,\vsi}
w_{\alpha} \exp \bigl( H_{N, t}(\vsi,\alpha) + H_N^{\mathrm{pert}}(
\vsi) \bigr). \label{Gint}
\end{equation}
If we define the Gibbs measure $\Gamma_t$ on $\Sigma_N\times\mathbb
{N}^{k}$ by
\[
\Gamma_t \bigl\{(\vsi,\alpha) \bigr\} \sim w_{\alpha}\exp
\bigl(H_{N,t}(\vsi,\alpha)+ H_N^{\mathrm{pert}}(\vsi) \bigr),
\]
then a straightforward calculation using Gaussian integration by parts gives
%
%
\begin{eqnarray}\label{phider2}
\varphi'(t) &=& -\tfrac{1}{2} \theta(1) +\tfrac{1}{2} \e
\bigl\langle \theta(q_{\alpha^1\wedge\alpha^2}) \bigr\rangle_{\Gamma_t}
\nonumber
\\[-8pt]
\\[-8pt]
\nonumber
&&{} -
\tfrac{1}{2} \e\bigl\langle \xi(R_{1,2})-R_{1,2}
\xi'(q_{\alpha^1\wedge\alpha^2}) +\theta (q_{\alpha
^1\wedge\alpha^2}) \bigr
\rangle_{\Gamma_t},
\end{eqnarray}
where $\langle \cdot\rangle_{\Gamma_t}$ is the Gibbs average with
respect to
$\Gamma_t^{\otimes2}.$ When $\xi(x)=\sum_{p\geq1}\beta_p^2 x^p$ does
not contain terms for odd $p\geq3$, $\xi$ is convex on $[-1,1]$, which
implies that the last term in (\ref{phider2}) is negative, and dropping
this term and integrating the corresponding inequality for $0\leq t\leq
1$, we obtain an upper bound on the free energy in~(\ref{Parisi}). The
fact that the representation of this upper bound in terms of the Ruelle
probability cascades coincides with the formula in (\ref{Pk}) is well
known and is explained in great detail in Chapter 14 in \cite{SG2}. If
the terms for odd $p\geq3$ are present, the function $\xi$ is only
convex on $[0,1]$, but the argument still works if we know that
$R_{1,2}$ is nonnegative with high probability under $\e\Gamma_t^{\otimes2}$.
This is where the perturbation term in (\ref{Gint})
comes into play to ensure that the Ghirlanda--Guerra identities hold
along the interpolation and, as a consequence, to ensure the positivity
of the overlap via Talagrand's positivity principle (see Section 12.3
in \cite{SG2}). In fact, an observation in \cite{Posit} shows that the
perturbation term $H_{N}^{\mathrm{pert}}$ forces the positivity of the
overlap uniformly over all measures on $\Sigma_N$ in the following
sense. If given a measure $\nu_N$ on $\Sigma_N$ we define a random
probability measure $\hat{\nu}_N$ on $\Sigma_N$ by the change of
density $d\hat{\nu}_N(\vsi) \sim\exp H_N^{\mathrm{pert}}(\vsi) \,
d\nu_N(\vsi
)$, then Theorem 1 in \cite{Posit} implies that for any $\eps>0$,
%
%
\begin{equation}
\lim_{N\to\infty}\sup_{\nu_N} \e\hat{\nu}_N^{\otimes2}
(R_{1,2} \leq-\eps) = 0.
\end{equation}
Using this for the marginal $\nu_N$ on $\Sigma_N$ of the Gibbs measure
$\gamma_t \{(\vsi,\alpha) \} \sim w_{\alpha}\exp H_{N,
t}(\vsi
,\alpha)$ on $\Sigma_N\times\Natural^k$ implies that the remainder
term in (\ref{phider2}) is asymptotically nonnegative and we can
proceed as in the case of even $p$-spin interactions.

\textit{The Aizenman--Sims--Starr scheme}. The proof of the lower bound
is done in several steps, but it begins with the Aizenman--Sims--Starr
scheme \cite{AS2}. Let us consider the Hamiltonian $H_N^{-}(\vsi) =
\sum_{p\geq1} \beta_p H_{N,p}^{-}(\vsi)$, where
%
%
\begin{equation}
H_{N,p}^{-}(\vsi) = \frac{1}{(N+1)^{(p-1)/2}} \sum
_{1\leq i_1,\ldots,i_p\leq N}g_{i_1,\ldots,i_p} \sigma_{i_1}\cdots
\sigma_{i_p}. \label{HamSKp2}
\end{equation}
Let $G_N^{-}$ and $\langle \cdot\rangle_{\mathunderscore}$ denote
the Gibbs measure and its average corresponding to the Hamiltonian
$H_N^{-} + H_N^{\mathrm{pert}}$ and let $z(\vsi)$ and $y(\vsi)$ be two
Gaussian processes on $\Sigma_N$ with covariances
%
%
\begin{equation}
\e z \bigl(\vsi^1 \bigr) z \bigl(\vsi^2 \bigr) =
\xi'(R_{1,2}),\qquad \e y \bigl(\vsi^1 \bigr) y
\bigl( \vsi^2 \bigr) = \theta(R_{1,2}) \label{covzy}
\end{equation}
independent of each other and all other random variables. Then the
Aizenman--Sims--Starr scheme in \cite{AS2} yields the following (see,
e.g., Section 15.8 in \cite{SG2}).

%
\begin{proposition}\label{PropAS2}
The lower limit of the free energy in (\ref{Parisi}) is bounded from
below by
%
%
\begin{equation}
\log2 + \liminf_{N\to\infty} \bigl( \e\log\bigl\langle \ch z(\vsi) \bigr
\rangle_{\mathunderscore} - \e\log\bigl\langle \exp y(\vsi) \bigr
\rangle_{\mathunderscore} \bigr). \label{AS2}
\end{equation}
\end{proposition}

The only difference here is that we included the perturbation term
$H_N^{\mathrm{pert}}(\vsi)$, but, since it is of a smaller order, one can
easily check that it does not affect the computation leading to this
representation. Below, we will express the limit~(\ref{AS2}) in terms
of some asymptotic Gibbs measure that satisfies the exact form of the
Ghirlanda--Guerra identities, but, in order to do so, we first need to
show that Propositions \ref{PropGG} and \ref{PropAS2} also hold with
nonrandom choices of the sequence $x=(x_p)_{p\geq1}$ (depending on
$N$) rather than on average over $x$. We mentioned above that the proof
of the Ghirlanda--Guerra identities is robust to modifications of the
Hamiltonian $H_N$ and, in particular, they hold for the Gibbs measure
$G_N^{-}$ so that if in (\ref{GG}) we replace $\langle \cdot\rangle
$ by $\langle
\cdot\rangle_{\mathunderscore}$, then (\ref{GGxlim})
still holds.
Let us consider a collection
\[
\FF= \bigl\{ (f,n,p) \dvtx p\geq1, n\geq2, f \mbox{ is a monomial of }
(R_{l,l'})_{l,l'\leq n} \bigr\}.
\]
Since this is a countable family, we can enumerate it,
$((f_j,n_j,p_j))_{j\geq1},$ and define a function
%
%
\begin{equation}
\phi_\FF= \phi_\FF(x) = \sum_{j\geq1}
2^{-j} \phi(f_j,n_j,p_j),
\label{FF}
\end{equation}
which depends on the variables in $x=(x_p)_{p\geq1}$. Since each
monomial $|f|\leq1$, we can see from the definition (\ref{GG}) that
$|\phi(f,n,p)|\leq2$ and, therefore, the Ghirlanda--Guerra identities
(\ref{GGxlim}) imply that $\e_x \phi_\FF\to0$. Let
%
%
\begin{equation}
\lambda= \lambda(x) = \e_g \log\bigl\langle \ch z(\vsi) \bigr
\rangle_{\mathunderscore} - \e_g \log\bigl\langle \exp y(\vsi) \bigr
\rangle_{\mathunderscore}, \label{lambda}
\end{equation}
where, again, $\e_g$ denotes the expectation with respect to all
Gaussian random variables for a fixed $x.$ We will need the following
simple lemma.
%
%
\begin{lemma} We can find $x=(x_p)_{p\geq1}$ such that
%
%
\begin{equation}
\phi_\FF(x) \leq2c (\e_x \phi_\FF)^{1/2}
\quad\mbox{and}\quad \lambda(x) \leq\e_x \lambda+ 2c (\e_x
\phi_\FF)^{1/2}, \label{both}
\end{equation}
where $c$ is a constant that depends only on the function $\xi.$
\end{lemma}

\begin{pf}
If we denote by $\e_z$ and $\e_y$ the expectations with
respect to $(z(\vsi))$ and $(y(\vsi))$, then (\ref{covzy}) and Jensen's
inequality imply
\[
0\leq\e_g \log\bigl\langle \ch z(\vsi) \bigr\rangle_{\mathunderscore}
\leq \e_g \log\bigl\langle \e_z \ch z(\vsi) \bigr
\rangle_{\mathunderscore} = \xi'(1)/2
\]
and
\[
0\leq\e_g \log\bigl\langle \exp y(\vsi) \bigr\rangle_{\mathunderscore}
\leq \e_g \log\bigl\langle \e_y \exp y(\vsi) \bigr
\rangle_{\mathunderscore} = \theta(1)/2
\]
and, therefore, $-c\leq\lambda(x) \leq c$ for $c=\xi'(1)+ \theta(1)$.
Given $\eps>0$, consider the event
\[
\Omega= \bigl\{x=(x_p)_{p\geq1} \dvtx\lambda(x) \leq
\e_x \lambda+\eps \bigr\}.
\]
Then, if $\mathbb{P}_x$ denotes the probability with respect to the i.i.d.
sequence $(x_p)_{p\geq1}$ with the uniform distribution on $[1,2]$,
\[
\e_x \lambda\geq(\e_x \lambda+\eps) \mathbb
{P}_x \bigl(\Omega^c \bigr) - c \mathbb{P}_x(
\Omega),
\]
and, therefore,
\[
\mathbb{P}_x(\Omega)\geq\frac{\eps}{\e_x \lambda+ \eps+c}>\frac
{\eps}{3c}
\]
for $\eps<c$. On the other hand, Chebyshev's inequality implies
\[
\mathbb{P}_x (\phi_\FF\leq\eps)\geq1-
\frac{\e_x
\phi_\FF}{\eps},
\]
and $\Omega\cap\{\phi_\FF\leq\eps\}\not= \varnothing$ if $\eps
/3c > \e_x \phi_\FF/\eps$. Taking $\eps= 2(c\e_x \phi_\FF
)^{1/2}$ (which is
$<c$ for large $N$) implies that we can find $x$ that satisfies both
inequalities in (\ref{both}).
\end{pf}

For each $N$, let us choose $x^N = (x_p^{N})_{p\geq1}$ that satisfies
(\ref{both}) and, since $\e_x \phi_\FF\to0$, we get
%
%
\begin{equation}
\lim_{N\to\infty}\phi_\FF \bigl(x^N \bigr) = 0 \quad\mbox{and}\quad \liminf_{N\to\infty} \e_x \lambda\geq\liminf_{N\to\infty}
\lambda \bigl(x^N \bigr). \label{xNlim}
\end{equation}
Let us redefine the Hamiltonian $H_N^{\mathrm{pert}}$ and the Gibbs measure
$G_N^{-}$ by fixing parameters $x=x^N$ and, since the measure now
depends only on the Gaussian randomness, we will write $\e$ instead of
$\e_g.$ By (\ref{xNlim}), Proposition \ref{PropAS2} still holds for
this redefined measure $G_N^{-}$ and, recalling (\ref{FF}),
%
%
\begin{equation}
\e\bigl\langle f R_{1,n+1}^p \bigr\rangle_{\mathunderscore} -
\frac{1}{n} \e \langle f \rangle_{\mathunderscore} \e\bigl\langle
R_{1,2}^p \bigr\rangle_{\mathunderscore} - \frac{1}{n}
\sum_{l=2}^{n}\e\bigl\langle f
R_{1,l}^p\bigr\rangle_{\mathunderscore} \to0 \label{GGk}
\end{equation}
for all $p\geq1, n\geq2$ and all monomials $f$ of
$(R_{l,l'})_{l,l'\leq n}$.\vadjust{\goodbreak}

\textit{Asymptotic Gibbs' measures}. Next, we will define an asymptotic
analogue of the Gibbs measure and represent the limit (\ref{AS2}) in
terms of this measure. Let $(\vsi^l)_{l\geq1}$ be an i.i.d. sample
from $G_N^{-}$ and let $R^N = (R_{l,l'}^N)_{l,l'\geq1}$ be the
normalized Gram matrix, or matrix of overlaps, of this sample. Consider
a subsequence $(N_k)$ along which the limit in (\ref{AS2}) is achieved
(now with nonrandom parameters $x^N$) and the distribution of $R^N$
under ${\e G_N^{-}}^{\otimes\infty}$ converges in the sense of
convergence of finite dimensional distributions to the distribution of
some array $R^\infty$. For simplicity of notation, let us assume that
the sequence $(N_k)$ coincides with natural numbers. Under ${\e
G_N^{-}}^{\otimes\infty}$, the array $R^N$ is weakly exchangeable,
which means that
%
%
\begin{equation}
\bigl(R^N_{\pi(l),\pi(l')} \bigr) \stackrel{d} {=}
\bigl(R^N_{l,l'} \bigr)
\end{equation}
for any permutation $\pi$ of finitely many indices. Obviously, this
property will be preserved in the limit so that $R^\infty$ is a weakly
exchangeable symmetric nonnegative definite array and, following \cite
{DS}, we will call any such array a Gram-de Finetti array. The
Dovbysh--Sudakov representation \cite{DS} then guarantees that all such
arrays are generated by i.i.d. samples from random measures on a
separable Hilbert space (see \cite{PDS} for a detailed proof).
%
%
\begin{proposition} \label{ThDS}
If $(R_{l,l'})_{l,l'\geq1}$ is a Gram-de Finetti array such that
$R_{l,l} = 1$, then there exists a random measure $G$ on the unit ball
of a separable Hilbert space such that
%
%
\begin{equation}
(R_{l,l'} )_{l,l'\geq1} \stackrel{d} {=} \bigl(\vrho^l
\cdot\vrho^{l'} + \delta_{l,l'} \bigl(1-\bigl\|\vrho^l
\bigr\|^2 \bigr) \bigr)_{l,l'\geq1}, \label{DS}
\end{equation}
where $(\vrho^l)$ is an i.i.d. sample from $G$.
\end{proposition}

The importance of the Dovbysh--Sudakov representation in spin glass
models was first clearly demonstrated in \cite{AA}, and other examples
where this representation played an important role can be found in
\cite
{ACh,PGG,PGGsimple} and \cite{Tal-New}. Let $G$ be a
random measure generating the array $R^\infty$, let $(\vrho^l)$ be an
i.i.d. sample from $G$ and let $R_{l,l'} = \vrho^l\cdot\vrho^{l'}$ for
$l\not= l'$ and $R_{l,l}=1$. For simplicity of notation, we will now
omit $\infty$ in~$R^\infty$. If we denote by $\langle \cdot\rangle
$ the average
with respect to $G$, then, by (\ref{GGk}), the measure $G$ satisfies
the Ghirlanda--Guerra identities,
%
%
\begin{equation}
\e\bigl\langle f R_{1,n+1}^p \bigr\rangle = \frac{1}{n}
\e\langle f \rangle \e\bigl\langle R_{1,2}^p\bigr\rangle +
\frac{1}{n} \sum_{l=2}^{n}\e\bigl
\langle f R_{1,l}^p \bigr\rangle \label{GGlimit}
\end{equation}
for all $p\geq1, n\geq2$ and all monomials $f$ of
$(R_{l,l'})_{l,l'\leq n}$. Approximating bounded functions of the
overlaps (in the $L^1$ sense) by polynomials, we also have
%
%
\begin{equation}
\e\bigl\langle f \psi(R_{1,n+1}) \bigr\rangle = \frac{1}{n} \e
\langle f \rangle \e\bigl\langle \psi(R_{1,2})\bigr\rangle +
\frac{1}{n} \sum_{l=2}^{n} \e\bigl
\langle f \psi(R_{1,l}) \bigr\rangle \label{GGgen}
\end{equation}
for bounded measurable functions $f$ and $\psi$. Below, the identities
(\ref{GGgen}) will allow us to identify these asymptotic Gibbs
measures, but, first, let us show how the limit in (\ref{AS2}) can be
represented in terms of $G$. By Theorem 2 in \cite{PGG}, (\ref{GGgen})
implies that if $q^*$ is the largest point in the support of the
distribution of $R_{1,2}$ under $\e G^{\otimes2}$, then $G$ is
concentrated on the sphere of radius $\sqrt{q^*}$ with probability one
and, therefore, $R$ is generated by $(\vrho^l\cdot\vrho^{l'} +
\delta_{l,l'}(1-q^*))_{l,l'\geq1}$. Let $z(\vrho)$ and $y(\vrho)$
be two
Gaussian processes on the unit ball of our Hilbert space with covariances
%
%
\begin{equation}
\e z \bigl(\vrho^1 \bigr) z \bigl(\vrho^2 \bigr) =
\xi' \bigl(\vrho^1\cdot\vrho^2 \bigr),\qquad \e y
\bigl(\vrho^1 \bigr) y \bigl(\vrho^2 \bigr) = \theta
\bigl( \vrho^1\cdot\vrho^2 \bigr), \label{zy}
\end{equation}
let $\eta$ be a standard Gaussian random variable independent of
everything else and let $\e_\eta$ denote the expectation in $\eta$
only. Then the following holds.
%
%
\begin{lemma}\label{ThG}
We have
%
%
\begin{equation}
\lim_{N\to\infty} \e\log\bigl\langle \ch z(\vsi)\bigr\rangle_{\mathunderscore} =
\e\log \e_\eta\bigl\langle \ch \bigl(z(\vrho) +\eta \bigl(
\xi'(1)- \xi' \bigl(q^* \bigr) \bigr)^{1/2}
\bigr) \bigr\rangle \label{limAS1}
\end{equation}
and
%
%
\begin{equation}
\qquad\lim_{N\to\infty} \e\log\bigl\langle \exp y(\vsi)\bigr\rangle_{\mathunderscore} =
\e\log \e_\eta\bigl\langle \exp \bigl(y(\vrho) +\eta \bigl(\theta(1)-
\theta \bigl(q^* \bigr) \bigr)^{1/2} \bigr)\bigr\rangle . \label{limAS2}
\end{equation}
\end{lemma}

The proof of Lemma \ref{ThG} is based on the following observation. For
a moment, let $R=(R_{l,l'})_{l,l'\geq1}$ be an arbitrary Gram-de
Finetti array such that $R_{l,l'} = 1$, let $\LL$ be its distribution
and let $G$ be any random measure generating $R$ as in Proposition \ref
{ThDS}. It is known that in some sense this measure is unique (see
Lemma 4 in \cite{PDS}), but we will not need it here. Let us define
%
%
\begin{equation}
\Phi(\LL) = \e\log\e_\eta\bigl\langle \ch \bigl(z(\vrho) +\eta \bigl(
\xi'(1)-\xi' \bigl(\|\vrho\|^2 \bigr)
\bigr)^{1/2} \bigr) \bigr\rangle . \label{Phi}
\end{equation}
The Gaussian process $z(\vrho)$ here is the same as in (\ref{zy}), but
we do not assume now that $G$ is concentrated on the sphere $\|\vrho\|^2 = q^*$. We will prove that the right-hand side in (\ref{Phi}) does
not depend on the choice of the measure $G$ and, indeed, depends only
on the distribution $\LL$ in a continuous fashion.
%
%
\begin{lemma}\label{LemPhi}
The function $\LL\to\Phi(\LL)$ defined in (\ref{Phi}) is continuous
with respect to weak convergence of the distribution $\LL$ .
\end{lemma}

\begin{pf*}{Proof of Lemma \ref{ThG}} Since $R^N$ is the Gram matrix of
the sequence $(N^{-1/2}\vsi^l)$, we can simply think of the measure
$G_N^{-}$ as defined on $N^{-1/2}\Sigma_N$ which is a subset of the
sphere $\| \vsi\|=1$ in $\Reals^N$. Then (\ref{covzy}) agrees with
(\ref
{zy}) and Lemma \ref{LemPhi} implies (\ref{limAS1}) since $R^N$
converges in distribution to $R^\infty$ and, as we mentioned above, the
Ghirlanda--Guerra identities (\ref{GGlimit}) imply that $G$ is
concentrated on the sphere $\|\vrho\|^2=q^*$. Equation (\ref{limAS2})
can be proved similarly.
\end{pf*}

\begin{pf*}{Proof of Lemma \ref{LemPhi}}
The proof is almost identical to the proof of Lemma 11 in \cite
{Pspins}. For simplicity of notation, let us denote
\[
z_\eta(\vrho) = z(\vrho) +\eta \bigl(\xi'(1)-
\xi' \bigl(\|\vrho\|^2 \bigr) \bigr)^{1/2}
\]
and let $\e_z$ be the expectation in the randomness of $(z(\vrho))$
conditionally on all other random variables. By standard concentration
inequalities for Gaussian processes (see, e.g., Lemma 3 in \cite
{Posit}), we have that for $a\geq1$,
%
%
\begin{equation}
\mathbb{P}_z \bigl(\bigl | \log\e_\eta\bigl\langle \ch
z_\eta( \vrho)\bigr\rangle - \e_z \log\e_\eta
\bigl\langle \ch z_\eta(\vrho) \bigr\rangle \bigr| \geq a \bigr)\leq\exp
\bigl(-ca^2 \bigr) \label{GaussC}
\end{equation}
for some small enough constant $c$ that depends only on the function
$\xi$ through~(\ref{zy}). Since
\[
0 \leq\e_z \log\e_\eta\bigl\langle \ch z_\eta(
\vrho)\bigr\rangle \leq\log\bigl\langle \e_z \e_\eta\ch
z_\eta(\vrho)\bigr\rangle = \xi'(1)/2,
\]
the inequality (\ref{GaussC}) implies that $\mathbb{P}( |\log\e_\eta\langle \ch
z_\eta(\vrho)\rangle | \geq a)\leq\exp(-ca^2)$ for small~$c$ and large
enough $a$ and, therefore, if we denote $\log_a x =\break  \max(-a, \min
(\log
x, a))$, then for large $a,$
%
%
\begin{equation}
\bigl| \e\log\e_\eta\bigl\langle \ch z_\eta(\vrho)\bigr\rangle -
\e \log_a \e_\eta\bigl\langle \ch z_\eta(\vrho)
\bigr\rangle \bigr| \leq\exp \bigl(-ca^2 \bigr). \label{partA}
\end{equation}
Next, if we define $\ch_a x = \min(\ch x,\ch a)$, then
using that $|\log_a x - \log_a y | \leq e^a |x-y|$ and $| \ch x - \ch_a x | \leq\ch x I(|x| \geq a)$, we can write
\begin{eqnarray*}
\bigl| \e\log_a \e_\eta\bigl\langle \ch z_\eta(\vrho)
\bigr\rangle - \e \log_a \e_\eta\bigl\langle
\ch_a z_\eta(\vrho) \bigr\rangle\bigr | &\leq& e^a \e
\bigl\langle \bigl| \ch z_\eta(\vrho) - \ch_a z_\eta(
\vrho) \bigr| \bigr\rangle
\\
&\leq& e^a \e\bigl\langle \ch z_\eta(\vrho) I
\bigl(\bigl|z_\eta( \vrho)\bigr| \geq a \bigr) \bigr\rangle .
\end{eqnarray*}
By H\"older's inequality, this can be bounded by
\[
e^a \bigl( \e\bigl\langle \e_{z,\eta} \ch^2
z_\eta(\vrho) \bigr\rangle \bigr)^{1/2} \bigl( \e\bigl\langle
\mathbb{P}_{z,\eta} \bigl(\bigl|z_\eta(\vrho)\bigr| \geq a \bigr) \bigr
\rangle \bigr)^{1/2} \leq\exp \bigl(-ca^2 \bigr)
\]
for small $c$ and large enough $a$ since $\mathbb{P}_{z,\eta}
(|z_\eta(\vrho)|
\geq a)\leq\exp(-ca^2).$
Combining with (\ref{partA}),
%
%
\begin{equation}
\bigl| \e\log\e_\eta\bigl\langle \ch z_\eta(\vrho)\bigr\rangle -
\e \log_a \e_\eta\bigl\langle \ch_a
z_\eta(\vrho) \bigr\rangle \bigr| \leq\exp \bigl(-ca^2 \bigr).
\label{partA2}
\end{equation}
Approximating the logarithm by polynomials on the interval
$[e^{-a},e^a]$,\break  $\e\log_a \e_\eta\langle \ch_a z_\eta
(\vrho
)\rangle $ can be approximated by a linear combination of moments
%
%
\begin{equation}
\e \bigl(\e_\eta\bigl\langle \ch_a z_\eta(\vrho)
\bigr\rangle \bigr)^r = \e\biggl\langle \e_z
\e_\eta\prod_{l\leq r} \ch_a
\bigl(z_{\eta^l} \bigl( \vrho^l \bigr) \bigr) \biggr\rangle ,
\label{preF}
\end{equation}
where we used replicas and where
\[
z_{\eta^l} \bigl(\vrho^l \bigr) = z \bigl(\vrho^l
\bigr) + \eta^l \bigl(\xi'(1)-\xi' \bigl(
\bigl\| \vrho^l\bigr\|^2 \bigr) \bigr)^{1/2}
\]
and $(\eta^l)$ are i.i.d. standard Gaussian. Since the covariance of
the Gaussian sequence $(z_{\eta^l}(\vrho^l))$ is equal to
\[
\xi' \bigl(\vrho^l\cdot\vrho^{l'} \bigr) +
\delta_{l,l'} \bigl(\xi'(1)-\xi' \bigl(\bigl\|
\vrho^l\bigr\|^2 \bigr) \bigr) = \xi'(R_{l,l'}),
\]
the function inside the Gibbs average on the right-hand side of (\ref
{preF}) is equal to
%
%
\begin{equation}
\e_z \e_\eta\prod_{l\leq r}
\ch_a \bigl(z_{\eta^l} \bigl(\vrho^l \bigr) \bigr) =
F \bigl((R_{l,l'})_{l,l'\leq r} \bigr) \label{Fdef}
\end{equation}
for some continuous bounded function $F$ of the overlaps
$(R_{l,l'})_{l,l'\leq r}$. Together with~(\ref{partA2}) this shows that
we can approximate $\Phi(\LL)$ arbitrarily well by a linear combination
of $\e\langle F(R)\rangle $ for some continuous bounded functions
$F$ of
finitely many overlaps, which proves that $\Phi(\LL)$ is continuous
with respect to the distribution $\LL$ of the overlap array $R$.
\end{pf*}

\textit{Identifying asymptotic Gibbs' measures using ultrametricity}.
To show that the lower bound in (\ref{AS2}) matches Guerra's upper
bound, it remains to identify the difference of (\ref{limAS1}) and
(\ref
{limAS2}) with the second and third terms of the functional~(\ref{Pk}).
Since the asymptotic Gibbs measure $G$ satisfies the Ghirlanda--Guerra
identities~(\ref{GGgen}), the main result in \cite{PU} implies that the
support of $G$ is ultrametric with probability one, that is,
%
%
\begin{equation}
\e\bigl\langle I \bigl(R_{1,2}\geq\min(R_{1,3},R_{2,3})
\bigr) \bigr\rangle =1. \label{ultra}
\end{equation}
Given $r\geq1$, let us consider a function $\kappa(q)$ on $[0,1]$
such that
%
%
\begin{equation}
\kappa(q) = j/r \qquad \mbox{for } j/r\leq q<(j+1)/r, j=0,\ldots, r-1
\label{finite}
\end{equation}
and $\kappa(1)=1$. Equation (\ref{ultra}) implies that for any $q$ the
inequality $q\leq\vrho^l\cdot\vrho^{l'}$ defines an equivalence
relation $l\sim l'$ and, therefore, the array $(I(q\leq
R_{l,l'}))_{l,l'\geq1}$ is nonnegative definite, since it is
block-diagonal with blocks consisting of all elements equal to one.
This implies that $R^\kappa= (\kappa(R_{l,l'}))_{l,l'\geq1}$ is
nonnegative definite since it can be written as a convex combination
\[
\kappa(R_{l,l'}) = \sum_{j=1}^{r}
\frac{1}{r}I \biggl(\frac
{j}{r}\leq R_{l,l'} \biggr).
\]
In addition, it is clear that $R^\kappa$ is weakly exchangeable and
satisfies the Ghirlanda--Guerra identities (\ref{GGgen}). Then, by the
Dovbysh--Sudakov representation (\ref{DS}), $R^\kappa$ can be generated
by a sample from some random measure $G^\kappa$ on the unit ball of a
Hilbert space. If for simplicity we assume that $q^*\not= j/r$ for
$j\leq r$, then $\kappa(q^*)$ is the largest point in the support of
the distribution of $\kappa(R_{1,2})$ and, by Theorem~2 in~\cite{PGG},
the measure $G^\kappa$ is concentrated on the sphere $\|\vrho\|^2 =
\kappa(q^*)$. When $r\to\infty$, the distribution of $R^\kappa$
converges weakly to the distribution of $R$ and if we denote by
$\langle
\cdot\rangle_\kappa$ the average with respect to the measure
$G^\kappa$,
then Lemma \ref{LemPhi} implies that
%
%
\begin{equation}
\e\log\e_\eta\bigl\langle \ch(z(\vrho) +\eta \bigl(
\xi'(1)- \xi' \bigl(\kappa \bigl(q^* \bigr)
\bigr)^{1/2} \bigr) \bigr\rangle_\kappa\label{fterm}
\end{equation}
approximates the right-hand side of (\ref{limAS1}). It is well known
that an ultrametric measure, such as $G^\kappa$,
that satisfies the Ghirlanda--Guerra identities and under which the
overlaps $R^\kappa$ take finitely many values as in (\ref{finite}), can
be identified with the discrete Ruelle probability cascades by the
Baffioni-Rosati theorem \cite{BR} (see the proof of Theorem 15.3.6 in
\cite{SG2} for details). The fact that in this case~(\ref{fterm})
coincides with $X_0(\vec{m},\vec{q})$ in (\ref{Pk}) with parameters
$k=r-1$, $q_j = j/r$ and $m_j = \e\langle I(R_{1,2}< q_{j+1})\rangle
$ is also
well known (see, e.g., Theorem~14.2.1 in~\cite{SG2}). One can similarly
show that~(\ref{limAS2}) corresponds to the second term in~(\ref{Pk})
and this finishes the proof of Theorem~\ref{Th1}. One could also work
with the continuous Ruelle probability cascades using a general theory
developed in~\cite{Bolthausen}, but, at this point, it was easier to
simply discretize the overlap array and Gibbs measure.

%



\printaddresses

\end{document}